\documentclass[a4paper,12pt,reqno]{amsart}

\usepackage{amssymb}

\edef\restoreparindent{\parindent=\the\parindent\relax}
\usepackage{parskip}
\restoreparindent

\makeatletter
\@namedef{subjclassname@2020}{\textup{2020} Mathematics Subject Classification}
\makeatother

\newtheorem{thm}{Theorem}
\newtheorem{lem}{Lemma}

\newtheorem{cor}{Corollary}
\newtheorem{conj}{Conjecture}

\newtheorem{Thm}{Theorem}

\theoremstyle{definition}
\newtheorem{rem}{Remark}
\newtheorem{defn}{Definition}

\newcommand{\IR}{{\mathbb R}}

\newcommand{\IC}{{\mathbb C}}
\newcommand{\ID}{{\mathbb D}}

\newcommand{\IT}{{\mathbb T}}
\newcommand{\real}{{\operatorname{Re}\,}}
\newcommand{\eit}{{e^{i\theta}}}

\newcommand{\be}{\begin{equation}}
\newcommand{\ee}{\end{equation}}

\newcommand{\blem}{\begin{lem}}
\newcommand{\elem}{\end{lem}}

\newcommand{\bdefn}{\begin{defn}}
\newcommand{\edefn}{\end{defn}}

\newcommand{\bthm}{\begin{thm}}
\newcommand{\ethm}{\end{thm}}

\newcommand{\bcor}{\begin{cor}}
\newcommand{\ecor}{\end{cor}}

\newcommand{\bconj}{\begin{conj}}
\newcommand{\econj}{\end{conj}}

\newcommand{\brem}{\begin{rem}}
\newcommand{\erem}{\end{rem}}

\newcommand{\bpf}{\begin{proof}}
\newcommand{\epf}{\end{proof}}

\newcommand{\ds}{\displaystyle}

\begin{document}
	
\bibliographystyle{abbrv}

\title[Gabriel's problem for harmonic Hardy spaces]
{Gabriel's problem for harmonic Hardy spaces}

\author{Suman Das}
\address{Department of Mathematics with Computer Science, Guangdong Technion - Israel
	Institute of Technology, Shantou, Guangdong 515063, P. R. China.}
\email{suman.das@gtiit.edu.cn}



\begin{abstract}
We obtain inequalities of the form $$\int_C |f(z)|^p\, |dz| \leq A(p) \int_{\IT}|f(z)|^p\, |dz|, \quad (p>1)$$ where $f$ is harmonic in the unit disk $\ID$, $\IT$ is the unit circle, and $C$ is any convex curve in $\ID$. Such inequalities were originally studied for analytic functions by R.~M.~Gabriel [Proc. London Math. Soc. 28(2), 1928]. We show that these results, unlike in the case of analytic functions, cannot be true in general for $0< p \le 1$. Therefore, we produce an inequality of a slightly different type, which deals with the case $0<p<1$. An example is given to show that this result is ``best possible", in the sense that an extension to $p=1$ fails. Then we consider the special case when $C$ is a circle, and prove a refined result which surprisingly holds for $p=1$ as well. We conclude with a maximal theorem which has potential applications.
\end{abstract}

\maketitle
\pagestyle{myheadings}
\markboth{Suman Das}{Gabriel's problem for harmonic Hardy spaces}


\section{Introduction and Main Result}\label{sec1}

Suppose $\ID=\{z \in {\mathbb C}:\, |z|<1\}$ is the open unit disk in the complex plane $\IC$, and $\mathbb{T}=\{z \in \IC : |z|=1\}$ is the unit circle. For a function $f$ analytic in $\ID$,
the integral mean $M_p(r, f)$ is defined as $$ M_p(r, f) = \left( \frac{1}{2\pi}\int_{0}^{2\pi}|f(r e^{i\theta})|^p\, d\theta \right)^{\frac{1}{p}}, ~~~ p>0.
$$
The function $f$ is said to be in the Hardy space $H^p$ if $M_p(r, f)$ remains
bounded as $r \to 1^-$. It is well known that every function $f \in H^p$ has a radial limit almost everywhere. We denote by $f(e^{i\theta})$ the radial limit of $f$ on $\IT$. For a detailed survey on Hardy spaces, we refer to the books of Duren \cite{Duren} and Koosis \cite{Koosis}.

In \cite{gab2}, R.~M.~Gabriel considered the integral of $|f|^p$ along certain curves, and proved the following result.

\begin{Thm}\label{theGab}{\rm \cite{gab2}}
	Let $f \in H^p$, $p >0$. If $C$ is any convex curve in $\ID$, then 
	$$
	\int_C |f(z)|^p \,|dz| \leq 2 \int_{\IT}|f(z)|^p \,|dz|.
	$$
	The constant $2$ on the right hand side is best possible.	
\end{Thm}

Indeed, the theorem was proved for an arbitrary circle $\Gamma$ (instead of $\IT$) and analyticity was assumed on and inside $\Gamma$. However, that hypothesis is readily replaced by the weaker condition $f\in H^p$.

This problem originated through a famous inequality of Riesz and Fej\'er, where $C$ is the diameter $-1\le x \le 1$ and the constant on the right hand side is $1/2$ (see Theorem 3.13 of \cite{Duren}). The study of similar inequalities continued in a series of individual papers by Gabriel and Frazer, including but possibly not limited to, \cite{gab1,gab2,gab3,gab4, fra1,fra3,fra5,fra6}. Later, Granados termed these as ``Gabriel's problem" in her review \cite{granados}, which contains the background and subsequent developments, as well as open problems.

In this paper, we study Gabriel's problem for complex-valued harmonic functions. A function $f=u+iv$ is harmonic in $\ID$, if $u$ and $v$ are real-valued harmonic functions in $\ID$. Every such function has a decomposition $f=h+\overline{g}$, where $h$ and $g$ are analytic functions in $\mathbb{D}$. This decomposition is unique up to an additive constant. We say that the function $f$ belongs to the harmonic Hardy space $h^p$ $(p>0)$ if $$\sup_{0<r<1} M_p(r,f) < \infty.$$

In recent times, this problem resurfaced when Kayumov, Ponnusamy and Sairam Kaliraj \cite{Kayumov_PonSai_1} obtained the harmonic analogue of the Riesz-Fej\'er inequality for the space $h^p$, $p \in (1,2]$.
\begin{Thm}\label{RF_har}\cite{Kayumov_PonSai_1}
If $f \in h^p$ for $1<p\le 2$, then
$$
\int_{-1}^{1}|f(x)|^p \,dx \leq \frac{1}{2} \sec^p\left( \frac{\pi}{2p} \right)\int_{0}^{2\pi}|f(e^{i\theta})|^p \,d\theta.
$$ The inequality is sharp.
\end{Thm}
Later Melentijevi\'{c} and Bo\v{z}in \cite{Melentijevic_Bozin} showed that this sharp result holds for $p>2$ as well. In \cite{DK1}, a slight variant of this inequality was established by the present author and Sairam Kaliraj. Very recently, Chen and Hamada \cite{ChenHam} have produced a Riesz-Fej\'er type inequality for pluriharmonic functions in $\IR^{2n}$, thereby extending the problem to higher dimensions.

Continuing in this direction, in this paper, we prove Gabriel's result (Theorem \ref{theGab}) for the harmonic Hardy space $h^p$. The main theorem is as follows.

\bthm\label{Dthm1} Let $f \in h^p$ for some $p>1$, and let $C$ be any convex curve in $\ID$. Then
$$ \int_{C} |f(z)|^p \,|dz| \le 4 \int_{\mathbb{T}} |f(z)|^p\, |dz|,
$$
if $p\ge 2$, and 
$$ \int_{C} |f(z)|^p\, |dz| \le 2\sec^p\left(\frac{\pi}{2p}\right) \int_{\mathbb{T}} |f(z)|^p \,|dz|,
$$ if $1<p<2$.
\ethm

We give the proof of this result, followed by its consequences, in Section \ref{sec2}. In Section \ref{sec3}, we consider the special case when $C$ is a circle, and produce a refinement of Theorem \ref{Dthm1}.

\section{Proof and Further Implications}\label{sec2}

\subsection{The case $p\ge 2$:} We start with the simple case $p=2$. If $f=h+\bar{g} \in h^2$, we observe that
\begin{align*}
\int_{\IT} |f(z)|^2\, |dz| = \int_{\IT} |h(z)|^2 \,|dz| + \int_{\IT} |g(z)|^2\, |dz|,
\end{align*}
so that $h,g \in H^2$. Therefore, Theorem \ref{theGab} gives
\begin{align*}
\int_{C} |f(z)|^2\, |dz| & \le \int_{C} (|h(z)|+|g(z)|)^2\, |dz|\\ & \le 2 \left[\int_{C} |h(z)|^2\, |dz| + \int_{C} |g(z)|^2\, |dz|\right] \\ & \le 4 \left[\int_{\IT} |h(z)|^2 \,|dz| + \int_{\IT} |g(z)|^2\, |dz|\right]\\ & = 4 \int_{\IT} |f(z)|^2 \,|dz|.
\end{align*}
Now assume that $f \in h^p$, $p>2$. Let us recall that $f$ has the Poisson integral representation
$$f(r\eit) = \frac{1}{2\pi} \int_{0}^{2\pi} \frac{1-r^2}{1-2r\cos(\theta-t)+r^2} f(e^{it}) \,dt.$$ Also, it is known that the Poisson integral of a function $\varphi \in L^1(\IT)$ is a harmonic function in $\ID$. For $z=r\eit$, an appeal to Jensen's inequality shows that
\begin{align*}
\int_{C} |f(z)|^p \,|dz| & \le \int_{C} \left(\frac{1}{2\pi} \int_{0}^{2\pi} \frac{1-r^2}{1-2r\cos(\theta-t)+r^2} |f(e^{it})| \,dt\right)^{\frac{p}{2}.2} \,|dz|\\ & \le \int_{C} \left(\frac{1}{2\pi} \int_{0}^{2\pi} \frac{1-r^2}{1-2r\cos(\theta-t)+r^2} |f(e^{it})|^{p/2} \,dt\right)^2\, |dz|.
\end{align*}
Let us write $$U(z) = \frac{1}{2\pi} \int_{0}^{2\pi} \frac{1-r^2}{1-2r\cos(\theta-t)+r^2} |f(e^{it})|^{p/2} \,dt, \quad z=r\eit.$$ Then $U$ is non-negative and harmonic in $\ID$ such that $U(\eit) = |f(\eit)|^{p/2}$. Therefore, from what we have already proved,
\begin{align*}
\int_{C} |f(z)|^p \,|dz| &\le \int_{C} U^2(z) \,|dz| \\ & \le 4 \int_{\IT} U^2(z)\, |dz|\\ & = 4 \int_{\IT} |f(z)|^p\, |dz|,
\end{align*}
as desired.

\subsection{The case $1<p<2$:} This requires little bit more technicalities. A positive real-valued function $u$ is called log-subharmonic, if $\log u$ is subharmonic. To prove this part of the theorem, we make use of the following classical result of Lozinski \cite{Lozinski}.

\begin{Thm}\label{Lozin} {\rm \cite{Lozinski}}
	Suppose $\Phi$ is a log-subharmonic function from $\ID$ to $\IR$, such that $$ \sup_{0<r<1} \int_{0}^{2\pi} \Phi^p(re^{i\theta})\, d\theta < \infty, \quad p>0.$$ Then there exists $f \in H^p$ such that $\Phi(z) \le |f(z)|$ for $z \in \ID$, and $\Phi(\eit) = |f(\eit)|$ for almost every $\theta$.
\end{Thm}
This allows us to prove the next lemma which is essential for our purpose.

\begin{lem}\label{Dlem}
	Let $f=h+\bar{g} \in h^p$ for some $p > 1$. Then
	$$\int_{C}(|h(z)|+|g(z)|)^p\, |dz| \leq 2 \int_{\mathbb{T}}(|h(z)|+|g(z)|)^p\, |dz|.$$
	The coefficient $2$ on the right-hand side is sharp.
\end{lem}
\bpf
If $f \in h^p$ for $p>1$, an easy argument shows that $h,g \in H^p$ (see, for example, Lemma 1 of \cite{DK2} for the more general result). This ensures that the integral on the right-hand side is bounded. Therefore, in order to apply Theorem \ref{Lozin}, it is enough to show that $\log(|h(z)| + |g(z)|)$ is subharmonic.

Let $z_0\in \ID$ be such that $h(z_0) \neq 0$. Let us write $$\log(|h(z)| + |g(z)|) = \log|h(z)| + \log\left(1+\left|\frac{g(z)}{h(z)}\right|\right)$$ in a neighbourhood of $z_0$. The function $\log|h(z)|$ is subharmonic. The second function is of the form $\log(1+e^G)$, where $G=\log|g/h|$ is subharmonic in a neighbourhood of $z_0$. The function $\varphi(x)=\log(1+e^x)$ is convex and increasing, this can be checked by taking the derivatives. It is well-known that the composition of a convex, increasing function with a subharmonic function is subharmonic. Therefore, so is $\varphi\circ G=\log(1+|g(z)/h(z)|)$. So, $\log(|h(z)| + |g(z)|)$ is the sum of two subharmonic functions, hence subharmonic.

Now, let $z_1\in \ID$ be a zero of $h$. If $g(z_1)\neq 0$, the above argument applies with $h$ and $g$ interchanged. If $g(z_1)=0$, then $\log(|h(z_1)| + |g(z_1)|)=-\infty$, so that the local sub-mean-value property readily holds at $z_1$. This shows that $\log(|h(z)| + |g(z)|)$ is subharmonic in the whole $\ID$.

Therefore, Theorem \ref{Lozin} implies that there exists $F \in H^p$ such that $$|h(z)|+|g(z)| \le |F(z)| \, \text{ in } \ID, \quad |h(\eit)|+|g(\eit)| = |F(\eit)| \, \text{ a.e.}$$ It follows that
\begin{align*}
\int_{C}(|h(z)|+|g(z)|)^p \,|dz| & \le \int_{C} |F(z)|^p \,|dz|\\ & \le 2 \int_{\IT} |F(z)|^p \,|dz| \quad (\text{Theorem \ref{theGab}})\\ & = 2 \int_{\mathbb{T}}(|h(z)|+|g(z)|)^p\, |dz|
\end{align*}
and the proof of the lemma is complete. The sharpness follows from Theorem \ref{theGab}.
\epf

We also require the following inequality of Kalaj, the proof of which employs a method of plurisubharmonic functions.
\begin{Thm}\label{thm_Kalaj_ana_har}{\rm \cite{Kalaj}}
	Let $1 < p < \infty$ and let $f = h +\overline{g} \in h^p$ with $\real(h(0)g(0))=0$. Then
	$$\int_{\mathbb{T}}(|h(z)|^2 +|g(z)|^2)^{p/2} \,|dz| \leq \frac{1}{\left(1 - |\cos \frac{\pi}{p}|\right)^{p/2}} \int_{\mathbb{T}}|f(z)|^p \,|dz| .$$
	The inequality is sharp.
\end{Thm}

Now we are ready to resume the proof of the case $1<p<2$. Suppose, without any loss of generality, $g(0)=0$. We see that
\begin{align*}
\int_{C}|f(z)|^p \,|dz| & \le \int_{C}(|h(z)|+|g(z)|)^p\, |dz|\\ & \le 2 \int_{\IT}(|h(z)|+|g(z)|)^p \,|dz| \quad (\text{Lemma \ref{Dlem}})\\ & \le 2^{\frac{p}{2}+1} \int_{\IT}(|h(z)|^2+|g(z)|^2)^{p/2}\, |dz|\\ & \le \frac{2^{\frac{p}{2}+1}}{\left(1 - |\cos \frac{\pi}{p}|\right)^{p/2}} \int_{\IT}|f(z)|^p \,|dz|,
\end{align*}
where the last inequality follows by Theorem \ref{thm_Kalaj_ana_har}. Finally, we see that, for $1<p<2$, $$\frac{2^{\frac{p}{2}+1}}{\left(1 - |\cos \frac{\pi}{p}|\right)^{p/2}} = 2 \sec^p \left(\frac{\pi}{2p}\right),$$ which completes the proof of the theorem.

\subsection{Sharpness and further ramifications}
This constant is best possible, in the sense that one cannot get rid of the factor $\sec^p (\frac{\pi}{2p})$. This can be easily seen from Theorem \ref{RF_har}, which is contained in Theorem \ref{Dthm1} (i.e., $C$ to be taken as the diameter $-1\le x \le 1$). Let us mention that the function $$f(z) = \real\left\{(1-z^2)^{-\frac{1}{p}}\right\} \quad (z \in \ID)$$ works as an extremal function in this case (see \cite{Melentijevic_Bozin}).

Also, we observe that $$\sec^p \left(\frac{\pi}{2p}\right) \to +\infty \ \textit{ as } \ p \to 1.$$ This allows us to show that, unlike in the case of analytic functions, results of this type are not true in general for $0<p\le 1$.

\bthm \label{Dthm2}
The inequality in {\rm Theorem \ref{Dthm1}} does not generally hold for the case $0<p\le 1$.
\ethm

\bpf
Suppose, on the contrary, that there is $p_1 \in (0,1]$ such that the inequality $$\int_{C} |f(z)|^{p_1} \,|dz| \le K(p_1) \int_{\mathbb{T}} |f(z)|^{p_1} \,|dz|$$ holds, where $K(p_1)$ is a constant depending only on $p_1$. For $p>p_1$ and $z=r\eit$, we use Poisson integral and Jensen's inequality to deduce that
\begin{align*}
	\int_{C} |f(z)|^p \,|dz| & \le \int_{C} \left(\frac{1}{2\pi} \int_{0}^{2\pi} \frac{1-r^2}{1-2r\cos(\theta-t)+r^2} |f(e^{it})| \,dt\right)^{\frac{p}{p_1}\cdot {p_1}}\, |dz|\\ & \le \int_{C} \left(\frac{1}{2\pi} \int_{0}^{2\pi} \frac{1-r^2}{1-2r\cos(\theta-t)+r^2} |f(e^{it})|^{p/p_1} \,dt\right)^{p_1}\, |dz|\\ & \le K(p_1) \int_{\mathbb{T}} |f(z)|^{\frac{p}{p_1}\cdot {p_1}}\, |dz|\\ & = K(p_1) \int_{\mathbb{T}} |f(z)|^p\, |dz|.
\end{align*}
We can now choose $p=1+\epsilon$, where $\epsilon>0$ is sufficiently small, to arrive at a contradiction.
\epf

However, we produce a workaround in the form of the next result, which deals with the case $0<p<1$.
\bthm\label{new}
Let $f\in h^1$ and $C$ be any convex curve in $\ID$. Then, for $0<p<1$, $$\int_C |f(z)|^p \,|dz| \leq A(p) \left(\int_{\IT}|f(z)| \,|dz|\right)^p,$$ where $$A(p)=2(2\pi)^{1-p}\left[1+\sec\left(\frac{\pi p}{2}\right)\right].$$
\ethm

\bpf
Let $U(z)$ be the Poisson integral of $|f(e^{it})|$. Then $U(z)$ is a non-negative harmonic function in $\ID$ with $U(e^{it})=|f(e^{it})|$. Since $|f(z)|$ is subharmonic in $\ID$, it follows (see \cite[p. 7]{Duren}) that $$|f(z)| \le U(z), \quad z \in \ID.$$ Suppose $V(z)$ is the harmonic conjugate of $U$ such that $V(0)=0$. Let us write $F(z)=U(z)+iV(z)$. Then, we see that
\begin{align*}
\int_C |f(z)|^p \,|dz| \le \int_C U^p(z)\, |dz| & \le \int_C |F(z)|^p \,|dz|\\ & \le 2 \int_{\IT} |F(z)|^p\, |dz| \quad (\text{Theorem \ref{theGab}})\\ & = 2 \int_{\IT} \left[U^2(z)+V^2(z)\right]^{p/2} \,|dz|\\ & \le 2 \left[\int_{\IT}U^p(z)\, |dz| + \int_{\IT}|V(z)|^p\, |dz|\right].
\end{align*}
As $1/p>1$, Jensen's inequality implies $$\left(\frac{1}{2\pi} \int_{\IT}U^p(z) \,|dz|\right)^{1/p} \le \frac{1}{2\pi} \int_{\IT}U(z) \,|dz|,$$ so that $$\int_{\IT}U^p(z)\, |dz| \le (2\pi)^{1-p} \left(\int_{\IT}U(z) \,|dz| \right)^p.$$ Also, a well-known theorem of Kolmogorov (see \cite[p. 57]{Duren}) leads to the inequality $$\int_{\IT}|V(z)|^p\, |dz| \le (2\pi)^{1-p} \sec\left(\frac{\pi p}{2}\right)\left(\int_{\IT}U(z) \,|dz| \right)^p.$$ Combining these, we find that
\begin{align*}
\int_C |f(z)|^p \,|dz| & \le 2(2\pi)^{1-p}\left[1+\sec\left(\frac{\pi p}{2}\right)\right] \left(\int_{\IT}U(z) \,|dz| \right)^p\\ & = A(p) \left(\int_{\IT}|f(z)| \,|dz| \right)^p,
\end{align*}
which completes the proof.
\epf

We are unable to check the sharpness of this result. Nevertheless, we show that $A(p)$ cannot be refined to remove the term $\sec\left(\frac{\pi p}{2}\right)$. Let $$f(z)=\real\left[\left(\frac{1+z}{1-z}\right)^p\right],$$ and suppose $C$ is the diameter $-1\le x \le 1$. It is obvious that $f \in h^1$, since the function $$L(z)=\frac{1+z}{1-z} \in H^{p}$$ for every $p<1$. Now, $$\int_C |f(z)|^p\, |dz| = \int_{-1}^1 \left(\frac{1+x}{1-x}\right)^{p^2}\, dx.$$ This integral converges by the Riesz-Fej\'er inequality. On the other hand, \begin{align*}
\int_{\IT}|f(z)|\, |dz| & = \int_{0}^{2\pi} \left|\frac{1+e^{it}}{1-e^{it}}\right|^p \left|\cos\left(p \arg \frac{1+e^{it}}{1-e^{it}}\right)\right| \,dt\\ & = \cos\left(\frac{\pi p}{2}\right) \int_{0}^{2\pi} \left|\frac{1+e^{it}}{1-e^{it}}\right|^p \,dt.\end{align*} The last integral is finite since $L\in H^p$. It follows that $$ \frac{\ds \int_C |f(z)|^p \,|dz|}{\ds \int_{\IT}|f(z)|\, |dz|} = K \sec \left(\frac{\pi p}{2}\right),$$ where $K$ is an absolute constant. Thus, $A(p)$ always contains a factor of $\sec\left(\frac{\pi p}{2}\right)$.

\brem
This example further confirms that no general result of this type is true for $p=1$, i.e.,  inequalities of the form $$\int_C |f(z)| \,|dz| \le B \int_{\IT} |f(z)| \,|dz| \quad (B>0 \text{ constant})$$ are typically not possible.
\erem

\section{A Special Case}\label{sec3}

In \cite{fra6}, Frazer produced the ``nicest" result of this type: when $C$ is a circle. In this case, the sharp constant turns out to be $1$.

\begin{Thm}\label{Frazer}\cite{fra6}
If $f \in H^p$ for some $p>0$ and $C$ is a circle in $\ID$, then
$$
\int_C |f(z)|^p \,|dz| \leq \int_{\IT}|f(z)|^p \,|dz|.
$$
Clearly the inequality is sharp.	
\end{Thm}

Here we obtain the harmonic analogue of this inequality, and show that the constant on the right hand side is always less than 2. Curiously, this result also holds for $p=1$. We do not know if this remains true for $0<p<1$, as the reasoning of Theorem \ref{Dthm2} fails.

\bthm\label{Dthm3}
Let $f \in h^p$ for some $p\ge 1$, and let $C$ be any circle in $\ID$. Then $$\int_{C} |f(z)|^{p}\, |dz| \le \int_{\mathbb{T}} |f(z)|^{p}\, |dz|,$$ if $p\ge2$, and $$\int_{C} |f(z)|^{p}\, |dz| < (1+r) \int_{\mathbb{T}} |f(z)|^{p} \,|dz|,$$ if $1\le p<2,$ where $r$ is the distance of the centre of $C$ from the origin. These inequalities are best possible.
\ethm

\bpf
The case $p \ge 2$ will be handled as before. If $f = h + \bar{g} \in h^2$, we see that
$$
\int_{C} |f(z)|^2\, |dz| =  \int_{C} |h(z)|^2\, |dz| + \int_{C} |g(z)|^2\, |dz| + 2 \real \int_{C} h(z)g(z) \,|dz|.
$$
The last integral is $0$ since $hg$ is analytic in $\ID$. Therefore,
\begin{align*}
	\int_{C} |f(z)|^2 \,|dz| & =  \int_{C} |h(z)|^2 \,|dz| + \int_{C} |g(z)|^2\, |dz| \\ & \le  \int_{\IT} |h(z)|^2 \,|dz| + \int_{\IT} |g(z)|^2\, |dz| \quad (\text{Theorem \ref{Frazer}})\\ & = \int_{\IT} |f(z)|^2\, |dz|.
\end{align*}
The result for $p>2$ can be easily obtained using Poisson integral and Jensen's inequality. The details are omitted.

The case $1\le p<2$ is more delicate. Suppose $z_0=r\eit$ is the centre of $C$ and the radius is $\rho < 1-|z_0|=1-r$. It is known that for $p\ge 1$, the function $|f(z)|^p$ is subharmonic in $\ID$. Let $U(z)$ be the Poisson integral of $|f(e^{it})|^p$. Then $U$ is non-negative and harmonic in $\ID$, with $U(e^{it})=|f(e^{it})|^p$. It follows, like in the proof of Theorem \ref{new}, that $$|f(z)|^p \le U(z), \quad z \in \ID.$$

The mean value property implies $$\int_{C} U(z) \,|dz| = 2\pi \rho U(z_0).$$ Therefore, using these and the Poisson integral representation of U, we get 
\begin{align*}
\int_{C}|f(z)|^p\, |dz| & \le \int_{C} U(z) \,|dz|\\ & = 2\pi \rho \left[\frac{1}{2\pi} \int_{0}^{2\pi} \frac{1-r^2}{1-2r\cos(\theta-t)+r^2} U(e^{it}) \,dt \right]\\ & = \rho \int_{0}^{2\pi} \frac{1-r^2}{1-2r\cos(\theta-t)+r^2} |f(e^{it})|^p \,dt \\ & < (1-r) \int_{0}^{2\pi} \frac{1-r^2}{(1-r)^2} |f(e^{it})|^p \,dt\\ & = (1+r) \int_{0}^{2\pi}  |f(e^{it})|^p \,dt\\ & = (1+r) \int_{\IT} |f(z)|^p\, |dz|.
\end{align*}
The sharpness can be seen by letting $r \to 0$, as $C$ and $\IT$ become concentric. This completes the proof.
\epf

We conclude this paper with another inequality which is closely related to the special case of circle. For this purpose, we recall the following result from \cite{DK1}, which is the consequence of a Riesz-Fej\'er type inequality.

\begin{Thm}\label{KDthm2}\cite{DK1}
Let $\{a_n\}$, $\{b_n\}$ be sequences of real numbers and $\theta$ be any acute angle. Then
$$
\sum_{k=0}^{\infty} \sum_{l=0}^{\infty} \frac{(a_ka_l+b_kb_l)\cos(k-l)\frac{\theta}{2}+2a_kb_l\cos(k+l)\frac{\theta}{2}}{k+l+1} \le \frac{2\pi}{\sin\frac{\theta}{2}+\cos\frac{\theta}{2}} \sum_{n=0}^{\infty}(a_n^2+b_n^2).$$
\end{Thm}
It is worthwhile to mention that for $a_n=b_n$ and $\theta=0$, this reduces to Hilbert's inequality. As an application of this, we obtain the next result which is somewhat along the line of, but far less intricate than, the Hardy-Littlewood maximal theorems.
\bthm\label{Dthm5}
If $f \in h^p$ for some $p\ge 2$, then
$$\int_{0}^1 \max_{|z|=r} |f(z)|^p\, dr \le \int_{\IT} |f(z)|^p\, |dz|.$$
\ethm

\bpf
It is enough to prove the theorem for $p=2$, as the the case $p>2$ can be deduced using our earlier techniques. Let us write $$f(z) = h(z)+\overline{g(z)} = \sum_{k=0}^{\infty} a_k z^k + \overline{\sum_{k=0}^{\infty} b_k z^k}.$$ It is easy to see that $$|f(r\eit)| \le \sum_{k=0}^{\infty} (|a_k|+|b_k|) r^k,$$ and therefore
\begin{align*}
\int_{0}^1 \max_{|z|=r} |f(z)|^2 \, dr& \le \int_{0}^1 \left[\sum_{k=0}^{\infty} (|a_k|+|b_k|) r^k\right]^2\, dr\\ & = \int_{0}^1 \left[\sum_{k=0}^{\infty} \sum_{l=0}^{\infty} (|a_k| |a_l| + |b_k| |b_l| + 2 |a_k| |b_l|) r^{k+l}\right] \, dr\\ & = \sum_{k=0}^{\infty} \sum_{l=0}^{\infty} \frac{|a_k| |a_l| + |b_k| |b_l| + 2 |a_k| |b_l|}{k+l+1}.
\end{align*}
It follows from Theorem \ref{KDthm2}, with $\theta=0$, that
\begin{align*}
\sum_{k=0}^{\infty} \sum_{l=0}^{\infty} \frac{|a_k| |a_l| + |b_k| |b_l| + 2 |a_k| |b_l|}{k+l+1} \le 2\pi \sum_{n=0}^{\infty}(|a_n|^2+|b_n|^2) = \int_{\IT} |f(z)|^2\, |dz|
\end{align*}
and the proof is complete.
\epf

\brem
This result can be particularly useful in finding the lower bounds of the $h^p$-norms of functions, defined as $$\|f\|_p = \left(\frac{1}{2\pi} \int_{\IT} |f(z)|^p \, |dz|\right)^{1/p}.$$ This is often more challenging than finding the upper bounds. For example, in the case of univalent harmonic functions, several upper bounds are given in \cite{DK2} and \cite{DK3}. But no such lower bounds are available in the literature, to the best of our knowledge.
\erem

\subsection*{Acknowledgement} The research was partially supported by Natural Science Foundation of Guangdong Province (Grant no. 2024A1515010467), and Li Ka Shing Foundation STU-GTIIT Joint-Research Grant (2024LKSFG06). The author is indebted to Antti Rasila and Anbareeswaran Sairam Kaliraj for useful suggestions that improved the quality of the paper.


\bibliography{references}


\end{document}